\documentclass[12pt]{article}

\usepackage{amsmath}
\usepackage{amsthm,amssymb}
\usepackage{multirow}
\usepackage{slashbox}
\usepackage{mathrsfs}
\usepackage{float}
\usepackage{color}
\usepackage{graphicx}
\usepackage{subfigure}
\usepackage{tabu}
\usepackage{booktabs}
\usepackage{array}
\usepackage{caption}
\usepackage{multirow}
\usepackage{geometry}
\usepackage{verbatim}

\makeatletter
\@addtoreset{equation}{section}
\makeatother

\newtheorem{de}{Definition}[section]
\newtheorem{thm}{Theorem}[section]
\newtheorem{lem}{Lemma}[section]
\newtheorem*{pr}{Proof}
\newtheorem{ex}{Example}[section]

\newcommand{\A}{{\mathcal A}}
\newcommand{\B}{{\mathcal B}}
\newcommand{\C}{{\mathcal C}}
\newcommand{\D}{{\mathcal D}}
\newcommand{\I}{{\mathcal I}}
\newcommand{\Q}{{\mathcal Q}}
\newcommand{\U}{{\mathcal U}}
\newcommand{\s}{{\mathcal S}}
\newcommand{\V}{{\mathcal V}}

\newcommand{\X}{{\mathcal X}}
\newcommand{\LL}{{\mathcal L}}
\newcommand{\M}{{\mathcal M}}
\newcommand{\E}{{\mathcal E}}
\newcommand{\Z}{{\mathcal Z}}
\newcommand{\W}{{\mathcal W}}

\newcommand{\Aj}{{\hat {\mathcal A}}}
\newcommand{\Bj}{{\hat {\mathcal B}}}

\newcommand{\Uj}{{\hat {\mathcal U}}}

\newcommand{\nnl}{{n \times n \times \ell}}
\newcommand{\nontnt}{{n_1 \times n_2 \times n_3}}

\newcommand{\pql}{{p \times q \times \ell}}

\newcommand{\myny}{{m_1 \times n_1}}
\newcommand{\meny}{{m_2 \times n_1}}

\newcommand{\mymy}{{m_1 \times m_1}}
\newcommand{\meme}{{m_2 \times m_2}}

\newcommand{\nyny}{{n_1 \times n_1}}
\newcommand{\nene}{{n_2 \times n_2}}

\newcommand{\nn}{{n \times n}}

\newcommand{\RR}{{\mathbb R}}
\newcommand{\CC}{{\mathbb C}}

\newcommand{\fold}{{\rm {fold}}}
\newcommand{\unfold}{{\rm {unfold}}}
\newcommand{\bcirc}{{\rm {bcirc}}}
\newcommand{\diag}{{\rm {diag}}}
\newcommand{\fft}{{\rm {fft}}}
\newcommand{\ifft}{{\rm {ifft}}}
\newcommand{\rank}{{\rm {rank}}}
\newcommand{\true}{{\rm {true}}}

\newcommand{\T}{{\rm {T}}}

\newcommand{\Ajj}{{\hat {A}}}
\newcommand{\Bjj}{{\hat {B}}}
\newcommand{\Cjj}{{\hat {C}}}
\newcommand{\Djj}{{\hat {D}}}

\newcommand{\Qjj}{{\hat {Q}}}
\newcommand{\Ujj}{{\hat {U}}}
\newcommand{\sjj}{{\hat {S}}}
\newcommand{\Vjj}{{\hat {V}}}

\newcommand{\Xjj}{{\hat {X}}}
\newcommand{\Zjj}{{\hat {Z}}}

\title{CS decomposition and GSVD for tensors based on the T-product}
\author{
Yating Zhang\footnote{School of Mathematical Sciences, Ocean University of China, Qingdao 266100, China.
E-Mail: {\tt zhangyating@stu.ouc.edu.cn}},
Xiaoxia Guo\footnote{School of Mathematical Sciences, Ocean University of China, Qingdao 266100, China. E-Mail: {\tt guoxiaoxia@ouc.edu.cn}},
Pengpeng Xie\footnote{Corresponding author: School of Mathematical Sciences, Ocean University of China, Qingdao 266100, China.
E-Mail: {\tt xie@ouc.edu.cn}. The work of this author is supported in part by NSFC grant 11801534.},
Zhengbang Cao\footnote{School of Mathematical Sciences, Ocean University of China, Qingdao 266100, China.
E-Mail: {\tt caozhengbang@stu.ouc.edu.cn}}.
}
\date{}          

\begin{document}

\maketitle                
\vspace{-1.2 cm}
\begin{abstract}
This paper derives the CS decomposition for orthogonal tensors (T-CSD) and the generalized singular value decomposition for two tensors (T-GSVD) via the T-product. The structures of the two decompositions are analyzed in detail and are consistent with those for matrix cases. Then the corresponding algorithms are proposed respectively. Finally, T-GSVD can be used to give the explicit expression for the solution of tensor Tikhonov regularization. Numerical examples demonstrate the effectiveness of T-GSVD in solving image restoration problems.
\\ \hspace*{\fill} 
\end{abstract}

\section{Introduction}

\hskip 2em In practical applications, we often need to deal with large-scale data sets. As the higher-order generalization of matrix, tensor can reserve the intrinsic structural information of the data in higher-dimensional space. Third-order tensors are of great interest. Especially, Kilmer, Martin and Perrone \cite{Kilmer08} defined a new type of tensor multiplication which is called T-product. Based on the T-product, many linear algebra tools can be generalized from two-dimensional space to higher-dimensional space. T-product has been proved to be useful in many areas, such as image and signal processing \cite{Kilmer13,Soltani16,Tarzanagh18}, computer vision \cite{Xie18,Yin19}, denoising \cite{Zhang18}, and many more.

\hskip 2em Another contribution of \cite{Kilmer08} is that they presented a new way to extend the matrix singular value decomposition (SVD) to tensors (T-SVD). Braman \cite{Braman10} proposed the definition of eigenvalues and eigenvectors of the third-order tensors. In \cite{Kilmer13}, Kilmer et al. showed that tensors have many properties similar to matrices. For example, they gave the definition of inner product between tensors and orthogonal projectors. In addition, they proposed the algorithms of QR, power iteration, Golub-Kahan iterative bidiagonalization and CG for tensors. The Moore-Penrose inverse of tensors, T-Jordan canonical form, tensor function, generalized tensor function and T-product semidefiniteness were investigated in \cite{Jin17,Miao19,Lund20,Miao20,Zheng20} respectively. Tensor-Tensor Product Toolbox \cite{Lu18,Lu18(2)} makes it convenient to perform the tensor computations in MATLAB. Recently, Reichel and Ugwu \cite{Reichel21} proposed many methods based on the T-Arnoldi process to solve linear discrete ill-posed problems defined by third-order tensors and the T-product. In \cite{He21}, He et al. generalized the quotient singular value decomposition (QSVD) and product singular value decomposition (PSVD) from matrices to tensors under the T-product. It is showed that their decompositions can be applied to color image watermarking process. Different from QSVD and PSVD, the generalized singular value decomposition (GSVD) proposed in \cite{Van Loan76} has been used more frequently in Tikhonov regularization for solving ill-posed least squares problems. This GSVD is closely related to the CS decomposition \cite{Paige81}. This motivates us to generalize the GSVD in \cite{Van Loan76} and the CS decomposition from matrices to tensors.

\hskip 2em In this paper, we develop new decomposition theories of third-order tensors. Specifically, our paper focuses on two aspects. One is to extend the matrix CS decomposition and GSVD to tensors by the T-product (T-CSD and T-GSVD). The other one is on the applications of the decomposition T-GSVD. By exploiting the T-GSVD, we can derive explicit expression for the solution of tensor Tikhonov regularization, which is applied to image restoration problems.

\hskip 2em This paper is organized as follows. In section 2, we review basic definitions and notations. In sections 3 and 4, we construct the T-CSD and T-GSVD, and discuss in detail the structures of the two decompositions. Then the algorithms are presented. Section 4 also discusses how to use T-GSVD for solving the tensor Tikhonov regularization problem. Three numerical examples are performed in section 5. Section 6 gives concluding remarks.


\section{Preliminaries}
\hskip 2em In this section, we give some basic notations and results.

\subsection{Notation and indexing}
\hskip 2em Throughout this paper, we use lowercase letters $a,b,\dotsc$ for scalars, lowercase bold letters $\mathbf{a},\mathbf{b},\dotsc$ for vectors, capital letters $A,B,\dotsc$ for matrices, and calligraphic letters $\A,\B,\dotsc$ for tensors. A third-order tensor $\A\in\CC^{\nontnt}$ is a multidimensional array with its $\left(n_1,n_2,n_3\right)$th entry $a_{ijk}$. The $i$th horizontal, lateral and frontal slices are denoted by $\A^{(i)}\equiv\A(i,:,:)$, $\overrightarrow{\A_i}\equiv\A(:,i,:)$ and $A_i\equiv\A(:,:,i)$ respectively.
We also use $\A_{i,j}\equiv\A(i,j,:)\in\CC^{1 \times 1 \times n_3}$ to denote its $\left(i,j\right)$th tubal scalar. The command unfold and fold are defined as
\begin{equation*}
   \unfold\left(\A\right)=
   \begin{bmatrix}
       A_1    \\
       A_2    \\
       \vdots  \\
       A_{n_3}
   \end{bmatrix},\qquad \fold\left(\unfold\left(\A\right)\right)=\A.
\end{equation*}

\subsection{Discrete Fourier transform}
\hskip 2em It is well known that block circulant matrices can be block diagonalized
by the discrete Fourier transform (DFT). That is, suppose that $F_n=\left({\omega}^{\left(j-1\right)\left(k-1\right)}\right)\in\CC^{n \times n}$ is the DFT matrix, where $\omega=e^{-\frac{2\pi i}{n}}$ is an $n$-th root of unity with $i=\sqrt{-1}$, then
\begin{equation}
\left(F_{n_3} \otimes I_{n_1}\right)
\bcirc \left(\A\right)
\left(F_{n_3}^{-1} \otimes I_{n_2}\right)=\diag\left(\Ajj_1,\Ajj_2,\dotsc,\Ajj_{n_3}\right)
,
\label{block}
\end{equation}
where $\otimes$ denotes the Kronecker product and $\bcirc\left(\A\right)$ is a block circulant matrix of the form
\begin{equation*}
\bcirc\left(\A\right)=
\begin{bmatrix}
A_1          &A_{n_3}     &A_{n_3-1}   & \dotsc        & A_2 \\
A_2          &A_1         &A_{n_3}     & \dotsc        & A_3 \\
\vdots       &\ddots      &\ddots      & \ddots        & \vdots\\
A_{n_3}      &A_{n_3-1}   &\ddots      & A_2           & A_1
\end{bmatrix}.
\end{equation*}

\hskip 2em In essence, $\Aj\in\CC^{\nontnt}$ is a new tensor computed by taking fast Fourier transform (FFT) along each tubal scalar of $\A$, $i.e.$,
\begin{equation*}
\Aj\equiv\fft\left(\A,[\ ],3\right)=\fold
\begin{bmatrix}
     \Ajj_1    \\
     \Ajj_2    \\
     \vdots    \\
     \Ajj_{n_3}
\end{bmatrix}.
\end{equation*}
By (\ref{block}), it is not hard to show
\begin{equation}
\begin{cases}
  \Ajj_1    =A_1 +  A_2        + \dotsc +         A_{n_3}        \\
  \Ajj_2    =A_1 + \omega A_2        + \dotsc +{\omega}^{n_3-1}A_{n_3}        \\
  \quad\dotsc  \\
  \Ajj_{n_3}=A_1 + {\omega}^{n_3-1}A_2+ \dotsc +{\omega}^{(n_3-1)(n_3-1)}A_{n_3}
\end{cases}
\label{relation}
\end{equation}
Applying the inverse FFT along each tubal scalar of $\Aj$ gives $\A\equiv\ifft(\Aj,[\ ],3)$.

\subsection{Definitions and propositions}
\hskip 2em The following definitions and  propositions were introduced in \cite{Kilmer13,Kilmer08,Miao20}.
\begin{de} (T-product) Suppose $\A\in\CC^{\nontnt}$ and $\B\in\CC^{n_2 \times n \times n_3}$, then the T-product $\A*\B$ is the tensor in $\CC^{n_1 \times n \times n_3}$
\begin{equation*}
 \A*\B=\fold\left(\bcirc\left(\A\right)\cdot\unfold\left(\B\right)\right).
\end{equation*}
\end{de}

\begin{de} (Block Tensor) Suppose $\A\in\CC^{m_1 \times n_1 \times n}$, $\B\in\CC^{m_1 \times n_2 \times n}$, $\C\in\CC^{m_2 \times n_1 \times n}$ and $\D\in\CC^{m_2 \times n_2 \times n}$. The block tensor
\begin{equation*}
\begin{bmatrix}
 \A & \B \\
 \C & \D
\end{bmatrix}\in\CC^{(m_1+m_2) \times (n_1+n_2) \times n}
\end{equation*}
is defined by compositing the frontal slices of four tensors.
\end{de}

\hskip 2em By the definition of the T-product, matrix block and block matrix multiplication,
\cite{Miao20} gives the tensor block multiplication.

\begin{thm} Suppose $\A_1\in\CC^{m_1 \times n_1 \times n}$, $\B_1\in\CC^{m_1 \times n_2 \times n}$, $\C_1\in\CC^{m_2 \times n_1 \times n}$, $\D_1\in\CC^{m_2 \times n_2 \times n}$, $\A_2\in\CC^{n_1 \times r_1 \times n}$, $\B_2\in\CC^{n_1 \times r_2 \times n}$, $\C_2\in\CC^{n_2 \times r_1 \times n}$, $\D_2\in\CC^{n_2 \times r_2 \times n}$
are complex tensors, then we have
\begin{equation*}
\begin{bmatrix}
  \A_1 & \B_1\\
  \C_1 & \D_1
\end{bmatrix}*\begin{bmatrix}
  \A_2 & \B_2\\
  \C_2 & \D_2
\end{bmatrix}=\begin{bmatrix}
  \A_1*\A_2+\B_1*\C_2  & \A_1*\B_2+\B_1*\D_2  \\
  \C_1*\A_2+\D_1*\C_2  & \C_1*\B_2+\D_1*\D_2
\end{bmatrix}.
\end{equation*}
\end{thm}

\hskip 2em From the properties of FFT, one can easily get the following result about block tensor.
\begin{lem} Let $\A\in\RR^{m_1 \times n_1 \times n}$, $\B\in\RR^{m_1 \times n_2 \times n}$, $\C\in\RR^{m_2 \times n_1 \times n}$, $\D\in\RR^{m_2 \times n_2 \times n}$, then
\begin{equation*}
 \rm{fft}\left(
 \begin{bmatrix}
 \A & \B \\
 \C & \D
 \end{bmatrix}
 ,[\ ],3 \right)=
 \begin{bmatrix}
  \rm{fft}\left(\A,[\ ],3\right) &
  \rm{fft}\left(\B,[\ ],3\right)\\
  \rm{fft}\left(\C,[\ ],3\right) &
  \rm{fft}\left(\D,[\ ],3\right)
 \end{bmatrix}
 \in\CC^{(m_1+m_2) \times (n_1+n_2) \times n}.
\end{equation*}
\end{lem}

\begin{de} (Tensor Transpose) Suppose $\A\in\RR^{\nontnt}$, then $\A^{\T}$ is the $n_2\times n_1\times n_3$ tensor obtained by transposing each of the frontal slices and then reversing the order of transposed slices 2 through $n_3$.
\end{de}

\begin{lem} Suppose tensors $\A,\B$ and $\C$ are well-defined, then
\begin{equation}\C=\A+\B\iff\hat{C}_{i}=\hat{A}_i+\hat{B}_i,\label{sum}\end{equation}
\begin{equation}\C=\A*\B\iff\hat{C}_{i}=\hat{A}_i\hat{B}_i.\label{product}\end{equation}
\end{lem}
\hskip 2em Note that, (\ref{sum}) and (\ref{product}) show the relationship between the tensor operations and matrix operations, which play a core role in proving the theorems in this paper.

\begin{de} (Identity Tensor) The $\nnl$ identity tensor $\I_{nn\ell}$ is the tensor whose first frontal slice is the $\nn$ identity matrix, and whose other frontal slices are all zeros.
\end{de}

\begin{de} (Tensor Inverse) The tensor $\A\in\RR^{\nnl}$ has an inverse $\B$ provided that
\begin{equation*}
\A*\B=\B*\A=\I.
\end{equation*}
\end{de}

\begin{de} (Orthogonal Tensor) The tensor $\Q\in\RR^{\nnl}$ is orthogonal if \begin{equation*}\Q^{\T}*\Q=\Q*\Q^{\T}=\I.\end{equation*}
The tensor $\Q\in\RR^{\pql}$ is partially orthogonal if \begin{equation*}\Q^{\T}*\Q=\I_{qq\ell}.\end{equation*}
\end{de}

\begin{de} (F-diagonal Tensor) We say a tensor is f-diagonal if its each frontal slice is diagonal.
\end{de}

\begin{thm} (T-SVD) Suppose $\A\in\RR^{\nontnt}$. Then it can be factorized as
\begin{equation*}
\A=\U*\s*\V^{\T},
\end{equation*}
where $\U\in\RR^{n_1\times n_1\times n_3}$, $\V\in\RR^{n_2\times n_2\times n_3}$ are orthogonal, and $\s\in\RR^{\nontnt}$ is an f-diagonal tensor.
\end{thm}

\section{CS decomposition for orthogonal tensors}
\hskip 2em In this section, we consider the CS decomposition for orthogonal tensors based on the T-product (T-CSD). We first derive a special case of T-CSD.

\begin{thm} (T-CSD (Thin version)) Consider the tensor
\begin{equation*}
\Q=\begin{bmatrix}
 \Q_1\\
 \Q_2
\end{bmatrix}
,\qquad\Q_1\in\RR^{\myny \times n},\Q_2\in\RR^{\meny \times n},
\end{equation*}
where $m_1\geq n_1$ and $m_2\geq n_1$. If $\Q$ is partially orthogonal, then there exist orthogonal tensors $\U\in\RR^{\mymy \times n}$, $\V\in\RR^{\meme \times n}$ and $\Z\in\RR^{\nyny \times n}$ such that
\begin{equation}
\begin{bmatrix}
   \U      & \mathcal{O}        \\
    \mathcal{O}      & \V
\end{bmatrix}^{\T}*
  \begin{bmatrix}
      \Q_1\\
      \Q_2
  \end{bmatrix}
  *\Z=
         \begin{bmatrix}
            \C  \\
            \s
         \end{bmatrix},
         \label{CS}
\end{equation}
where $\C, \s$ are f-diagonal tensors, and satisfy
\begin{equation}
\C^{\T}*\C+\s^{\T}*\s=\I.
\label{CSI}
\end{equation}
\end{thm}

\begin{pr}
\upshape{Since} $\Q$ is partially orthogonal, $\Q^{\T}*\Q=\I$ follows from Definition 2.6. We then take FFT to both sides of this equation. From (\ref{product}), we have $\Qjj_i^*\Qjj_i=I_{n_1}$, $i.e.$ each $\Qjj_i$ is an orthonormal matrix, $i=1,2,\dotsc,n$. By Lemma 2.1, we partition the matrix
\begin{equation*}
\Qjj_i=
\begin{bmatrix}
 \left(\Qjj_1\right)_i \\
 \left(\Qjj_2\right)_i
\end{bmatrix},\ \left(\Qjj_1\right)_i\in\CC^{m_1 \times n_1},\ \left(\Qjj_2\right)_i\in\CC^{m_2 \times n_1},
\end{equation*}
where $m_1 \geq n_1$ and $m_2 \geq n_1$. Now, we can compute the CS decomposition \cite[Theorem 2.5.2]{Golub13} of matrix $\Qjj_i$. Thus, there exist unitary matrices $\Ujj_i\in\CC^{m_1 \times m_1}$, $\Vjj_i\in\CC^{m_2 \times m_2}$ and $\Zjj_i\in\CC^{n_1 \times n_1}$ such that
\begin{equation}
\begin{bmatrix}
\Ujj_{i}   & 0           \\
   0       & \Vjj_{i}
\end{bmatrix}^{*}
  \begin{bmatrix}
  \left(\Qjj_{1}\right)_i\\
  \left(\Qjj_{2}\right)_i
  \end{bmatrix}\Zjj_{i}=
      \begin{bmatrix}
      \Cjj_i\\
      \sjj_i
      \end{bmatrix},i=1,2,\dotsc,n,
      \label{CS_matrix}
\end{equation}
where $\Cjj_i$ and $\sjj_i$ are diagonal matrices of the following form
\begin{equation}
\begin{split}
\Cjj_i=\rm{diag}\left(\cos{\theta_{i,1}},\cos{\theta_{i,2}},\dotsc,\cos{\theta_{i,n_1}}\right)
      \in\RR^{\myny},\\
\sjj_i=\rm{diag}\left(\sin{\theta_{i,1}},\sin{\theta_{i,2}},\dotsc,\sin{\theta_{i,n_1}}\right)
      \in\RR^{\meny},
\end{split}
\label{C_S}
\end{equation}
and
\begin{equation}
\Cjj_i^{\T}\Cjj_i+\sjj_i^{\T}\sjj_i=I_{n_1}, i=1,2,\dotsc,n.
\label{CSI_matrix}
\end{equation}

If we define $\Uj$ as a tensor with its $i$th frontal slice being unitary matrix $\Ujj_i$, and $\U=\ifft(\Uj,[\ ],3)$, then $\U$ is an orthogonal tensor. Similarly, orthogonal tensors $\V$, $\Z$ and f-diagonal tensors $\C$, $\s$ can be constructed. Finally, by using (\ref{product}) again, the results (\ref{CS}) and (\ref{CSI}) follow from (\ref{CS_matrix}) and (\ref{CSI_matrix}) respectively.\qedsymbol
\end{pr}

\hskip 2em We summarize the process of the above proof in Algorithm 3.1.
\begin{table}[htb]
\centering  
\begin{tabular}{l}
\toprule
\textbf{Algorithm 3.1:} Compute the T-CSD for a partially orthogonal tensor $\Q$\\ \toprule
\quad\textbf{Input:} $\Q\in\RR^{\left(m_1+m_2\right) \times n_1 \times n}$\\
\quad\textbf{Output:} $\C\in\RR^{m_1\times n_1\times n}$, $\s\in\RR^{m_2 \times n_1 \times n}$,\\
\qquad\qquad\ \ \quad $\U\in\RR^{m_1 \times m_1 \times n}$, $\V\in\RR^{m_2 \times m_2 \times n}$, $\Z\in\RR^{n_1 \times n_1 \times n}$,\\
\qquad\qquad\ \ \quad such that $\Q_1=\U*\C*\Z^{\T}$, $\Q_2=\V*\s*\Z^{\T}$\\ \\
\quad$\Q=\fft\left(\Q,[\ ],3\right)$\\
\quad\textbf{for} $i=1$ \textbf{to} $n$\\
\quad\qquad $\left[\U\left(:,:,i\right),\V\left(:,:,i\right),\Z\left(:,:,i\right),\C\left(:,:,i\right),\s\left(:,:,i\right)\right]$\\
\quad\qquad\qquad\qquad\qquad\qquad\qquad$=\rm{csd}\left(\Q\left(1:m_1,:,i\right),\Q\left(m_1+1:m_1+m_2,:,i\right)\right)$\\
\quad\textbf{end for}\\
\quad$\U=\ifft\left(\U,[\ ],3\right)$, $\V=\ifft\left(\V,[\ ],3\right)$, $\Z=\ifft\left(\Z,[\ ],3\right)$, \\
\quad$\C=\ifft\left(\C,[\ ],3\right)$, $\s=\ifft\left(\s,[\ ],3\right)$\\
\toprule
\end{tabular}
\end{table}
\vspace{-0.7cm}

\hskip 2em By using matrix CS decomposition \cite[Theorem 2.5.3]{Golub13} and the same techniques in Theorem 3.1, it is possible to prove the more general version of the T-CSD.

\begin{thm} (T-CSD) Suppose
\begin{equation*}
\Q = \bordermatrix{ ~  &  n_1      & n_2   \cr
                    m_1& \Q_{11} & \Q_{12} \cr
                    m_2& \Q_{21} & \Q_{22} \cr }
\end{equation*}
is an orthogonal tensor in $\RR^{(m_1+m_2) \times (n_1+n_2) \times n}$ and that $m_1\geq n_1$ and $m_1\geq m_2$. Define the nonnegative integers $p$ and $q$ by $p=\max\{0,n_1-m_2\}$ and $q=\max\{0,m_2-n_1\}$. There exist orthogonal tensors $\U\in\RR^{\mymy \times n}$, $\V\in\RR^{\meme \times n}$, $\W\in\RR^{\nyny \times n}$ and $\Z\in\RR^{\nene \times n}$ such that
\begin{equation*}
\begin{bmatrix}
\U                &  \mathcal{O}\\
\mathcal{O}       &\V
\end{bmatrix}^{\T}*\Q*\begin{bmatrix}
\W                 &\mathcal{O}\\
\mathcal{O}        &\Z
\end{bmatrix}
=
\D
=
\bordermatrix{~   &  p        &   n_1-p    &   n_1-p    &    q        &  m_1-n_1  \cr
             \quad\, p  &  \I & \mathcal{O}&\mathcal{O} &  \mathcal{O}&\mathcal{O}\cr
                 n_1-p  &\mathcal{O}&\C    &     \s     &  \mathcal{O}&\mathcal{O}\cr
                 m_1-n_1&\mathcal{O}&\mathcal{O}&\mathcal{O}&\mathcal{O}   &\I    \cr
                 n_1-p  &\mathcal{O}&\s    &-\C       &\mathcal{O}   &\mathcal{O} \cr
              \quad\, q &\mathcal{O}&\mathcal{O}&\mathcal{O} &\I     &\mathcal{O} \cr   }
\end{equation*}
where $\C, \s$ are f-diagonal tensors and satisfy
\begin{equation*}
\C^{\T}*\C+\s^{\T}*\s=\I.
\end{equation*}
\end{thm}

\hskip 2em Note that, when $n=1$, Theorem 3.1 and Theorem 3.2 correspond to the matrix CS decomposition of thin version and general version respectively.


\section{GSVD for two tensors and its application}
\hskip 2em Recently in \cite{He21}, He et al. extended QSVD \cite{Moor89} and PSVD \cite{Heate86} from matrices to tensors. Here, QSVD and PSVD are two different generalizations of SVD, and they are now also known as two kinds of GSVD. Before that, Van Loan extended the traditional SVD for single matrix to propose a construction method for decomposing two matrices simultaneously for the first time in 1976 \cite{Van Loan76}. The decomposition proposed by Van Loan is commonly referred to as GSVD, which is one of the essential tools in numerical linear algebra. Naturally, we are interested in whether this GSVD can be generalized from matrices to tensors, and what kind of structure it will have.

\hskip 2em Meanwhile, Tikhonov regularization method \cite{Tikhonov63} is a well-known and highly regarded method for solving discrete ill-posed problems. GSVD is a useful tool for analysis of Tikhonov regularization problem. In \cite{Reichel21}, Reichel and Ugwu introduced linear ill-posed tensor least squares problems. Moreover, they transformed solving this kind of problems into solving one penalized least squares problem, which is called tensor Tikhonov regularization. Thus, in this section, we will establish the GSVD for two tensors via the T-product (T-GSVD). Then the role of T-GSVD in solving Tikhonov regularization problem is analyzed.

\subsection{T-GSVD}
\hskip 2em The detail of T-GSVD is as follows.
\begin{thm} (T-GSVD) Suppose $\A\in\RR^{m_1 \times n_1 \times n}$, $\B\in\RR^{m_2 \times n_1 \times n}$ with $m_1\geq n_1$, then there exist orthogonal tensors $\U\in\RR^{m_1 \times m_1 \times n}$ and $\V\in\RR^{m_2 \times m_2 \times n}$ and invertible tensor $\mathcal{X}\in\RR^{n_1 \times n_1 \times n}$ such that
\begin{equation*}
\U^{\T}*\A*\X=\D_{\A},\qquad\V^{\T}*\B*\X=\D_{\B},
\end{equation*}
where $\D_{\A}$ and $\D_{\B}$ are f-diagonal tensors.
\end{thm}

\begin{pr} \upshape{First, we compute} $\Aj$ and $\Bj$ by using (\ref{block}). Then, we can apply GSVD \cite[Theorem 2]{Van Loan76} to each matrix pair $\{\Ajj_i,\Bjj_i\}$. Since $\Ajj_i\in\CC^{m_1 \times n_1}$, $\Bjj_i\in\CC^{m_2 \times n_1}$ with $m_1\geq n_1$, there exist unitary matrices $\Ujj_i\in\CC^{m_1 \times m_1}$, $\Vjj_i\in\CC^{m_2 \times m_2}$ and invertible matrices $\Xjj_i\in\CC^{n_1 \times n_1}$, such that
\begin{equation}
\Ujj_i^{*}\Ajj_i\Xjj_i=\left(\Djj_{\A}\right)_i\!\!=\!\!\begin{bmatrix}
&\left(\sjj_{\A}\right)_i   & 0{\ }           & 0{\quad}\\
& 0                         & I_{p_i}{\ }     & 0{\quad} \\
& 0                         & 0{\ }           & 0{\quad} \\
\end{bmatrix},
\left(\sjj_{\A}\right)_i\!=\!\diag\left(\alpha_{p_i+1}^i,\dotsc,\alpha_{r_i}^{i}\right),
\label{GSVD1}
\end{equation}
\begin{equation}
\Vjj_i^{*}\Bjj_i\Xjj_i=\left(\Djj_{\B}\right)_i\!\!=\!\!\begin{bmatrix}
&\left(\sjj_{\B}\right)_i   & 0{\quad}           & 0{\quad} \\
& 0                         & 0{\quad}           & 0{\quad} \\
& 0                         & 0{\quad}           & 0{\quad} \\
\end{bmatrix},
\left(\sjj_{\B}\right)_i\!=\!\diag\left(\beta_{p_i+1}^i,\dotsc,\beta_{r_i}^{i}\right),
\label{GSVD2}
\end{equation}
where diagonal matrices $(\sjj_{\A})_i$ and $(\sjj_{\B})_i$ satisfy
\begin{equation}
\left(\sjj_{\A}\right)_i^{\T}\left(\sjj_{\A}\right)_i+\left(\sjj_{\B}\right)_i^{\T}\left(\sjj_{\B}\right)_i=I_{r_i-p_i},
\label{GSVD3}
\end{equation}
where
\begin{equation*}
p_i=\max\left\{r_i-m_2,0\right\}, r_i=\rank\left(
\begin{bmatrix}
\Ajj_i   \\
\Bjj_i
\end{bmatrix}
\right).
\end{equation*}
Now, we denote
\begin{equation*}
\Uj=\fold
\begin{bmatrix}
\Ujj_1\\
\Ujj_2\\
\vdots\\
\Ujj_n
\end{bmatrix},
\U=\ifft\left(\Uj,[\ ],3\right).
\end{equation*}
Obviously, $\U$ is an orthogonal tensor. Analogously, we can get orthogonal tensor $\V$, invertible tensor $\X$, and f-diagonal tensors $\s_{\A},\ \s_{\B},\ \D_{\A},\ \D_{\B}$.
By the inverse FFT, it is easy to show the frontal slices of $\D_{\A}$ and $\D_{\B}$ have the following form
\begin{gather*}
\left(D_{\A}\right)_{i} = \frac{1}{n}\sum^{n}_{j=1}\bar{\omega}^{(i-1)(j-1)}\left(\Djj_{\A}\right)_j,\\
\left(D_{\B}\right)_{i} = \frac{1}{n}\sum^{n}_{j=1}\bar{\omega}^{(i-1)(j-1)}\left(\Djj_{\B}\right)_j.
\end{gather*}
Finally, combining (\ref{product}), (\ref{GSVD1})-(\ref{GSVD3}) and notations defined above, we get the T-GSVD for the tensor pair $\left\{\A,\B\right\}$.\qedsymbol
\end{pr}

\hskip 2em Note that, when $n=1$, the results in Theorem 4.1 coincide with the GSVD for matrix cases.

\hskip 2em As we can see from (\ref{C_S}), the dimension of diagonal matrices $\Cjj_i$ and $\sjj_i$ relies only upon the dimension of the unitary matrix $\Qjj_i$. Once the orthogonal tensor $\Q$ is given, the dimension of matrix $\Qjj_i$ is fixed. Hence, the diagonal matrices $\Cjj_i$ and $\sjj_i$ are of the same dimension respectively. It follows that $\C^{\T}*\C+\s^{\T}*\s=\I$. Unfortunately, we cannot obtain the same result while computing T-GSVD. From (\ref{GSVD1}) and (\ref{GSVD2}), the dimension of $(\sjj_{\A})_i$ and $(\sjj_{\B})_i$ is not only related to the dimension but also the rank of the block matrix $\begin{bmatrix} \Ajj_i^{*}\ ,\Bjj_i^*\end{bmatrix}$. That is, the dimension of the identity matrix in the right side of (\ref{GSVD3}) cannot be guaranteed to be identical. Consequently, $\s_{\A}^{\T}*\s_{\A}+\s_{\B}^{\T}*\s_{\B}=\I$ does not hold directly. In other words, the equality holds if and only if $r_1=\cdots=r_n$. In many image restoration problems, this condition can be satisfied.

\hskip 2em The proof of Theorem 4.1 implies the specific algorithm for computing T-GSVD.

\begin{table}[htb]
\centering  
\begin{tabular}{l}
\toprule
\textbf{Algorithm 4.1:} Compute the T-GSVD for two tensors $\A$ and $\B$\\ \toprule
\quad\textbf{Input:} $\A\in\RR^{m_1 \times n_1 \times n}$, $\B\in\RR^{m_2 \times n_1 \times n}$\\
\quad\textbf{Output:} $\D_{\A}\in\RR^{m_1\times n_1\times n}$, $\D_{\B}\in\RR^{m_2 \times n_1 \times n}$,\\
\qquad\qquad\ \ \quad $\U\in\RR^{m_1 \times m_1 \times n}$, $\V\in\RR^{m_2 \times m_2 \times n}$, $\X\in\RR^{n_1 \times n_1 \times n}$,\\
\qquad\qquad\ \ \quad such that $\A=\U*\D_{\A}*\X^{-1}$, $\B=\V*\D_{\B}*\X^{-1}$\\ \\
\quad$\A=\fft\left(\A,[\ ],3\right)$, $\B=\fft\left(\B,[\ ],3\right)$\\
\quad\textbf{for} $i=1$ \textbf{to} $n$\\
\quad\qquad $\left[\U\left(:,:,i\right),\V\left(:,:,i\right),\X\left(:,:,i\right),\D_{\A}\left(:,:,i\right),\D_{\B}\left(:,:,i\right)\right]$\\
\quad\qquad\qquad\qquad\qquad\qquad\qquad\qquad\qquad\qquad$=\rm{gsvd}\left(\A\left(:,:,i\right),\B\left(:,:,i\right)\right)$\\
\quad\textbf{end for}\\
\quad$\U=\ifft\left(\U,[\ ],3\right)$, $\V=\ifft\left(\V,[\ ],3\right)$, $\X=\ifft\left(\X,[\ ],3\right)$, \\
\quad$\D_{\A}=\ifft\left(\D_{\A},[\ ],3\right)$, $\D_{\B}=\ifft\left(\D_{\B},[\ ],3\right)$\\
\toprule
\end{tabular}
\end{table}
\vspace{-0.3cm}

\subsection{Tensor Tikhonov regularization problem}
\hskip 2em Consider the linear ill-posed tensor least squares problem of the form
\begin{equation}
\min_{\overrightarrow{\X}\in\RR^{m \times 1 \times n}}\parallel\A*\overrightarrow{\X}-\overrightarrow{\B}\parallel_F,
\label{minF}
\end{equation}
where $\A\in\RR^{m \times m \times n}$ is a third-order tensor of ill-determined tubal rank \cite{Reichel21}, and $\parallel\cdot\parallel_F$ denotes the Frobenius norm of a third-order tensor
\begin{equation*}
\parallel\A\parallel_F=\sqrt{\sum_{i=1}^m\sum_{j=1}^m\sum_{k=1}^n a_{ijk}^2}.
\end{equation*}

\hskip 2em Reichel and Ugwu pointed out that (\ref{minF}) can be replaced by the penalized least squares problem
\begin{equation}
\min_{\overrightarrow{\X}\in\RR^{m \times 1 \times n}}\parallel\A*\overrightarrow{\X}-\overrightarrow{\B}\parallel_F^2+\mu^{-1}\parallel\LL*\overrightarrow{\X}\parallel_F^2,
\label{Tikhonov}
\end{equation}
where $\LL\in\RR^{s\times m\times n}$ is a regularization operator and $\mu>0$ is a regularization parameter. This replacement is called Tikhonov regularization. The two common forms of $\LL$ are $\LL_1\in\RR^{\left(m-2\right) \times m \times n}$ and $\LL_2\in\RR^{\left(m-1\right) \times m \times n}$, where
\begin{equation*}
\left(\LL_1\right)_1=\frac{1}{4}\begin{bmatrix}
-1      &2       &-1      &        &  \\
        &\ddots  &\ddots  &\ddots  &  \\
        &        &-1      &2       &-1\\
\end{bmatrix},\quad
\left(\LL_2\right)_1=\frac{1}{2}\begin{bmatrix}
1      &-1            &        &       &  \\
       &1             &-1      &       &  \\
       &              &\ddots  &\ddots &  \\
       &              &        &1      &-1\\
\end{bmatrix},
\end{equation*}
are the first frontal slices of $\LL_1$ and $\LL_2$ respectively. The other frontal slices are zero matrices. The methods on how to choose $\mu$ are introduced in detail in \cite{Reichel21}. Reichel and Ugwu proved theoretically that (\ref{Tikhonov}) has a unique solution
for any $\mu>0$ and they proposed many methods to solve (\ref{Tikhonov}). However, their methods are all iterative. With T-GSVD of the tensor pair $\left\{\A,\LL\right\}$, we can get an explicit formula for the unique solution.

\hskip 2em In case of $\LL=\LL_{1}\in\RR^{\left(m-2\right) \times m \times n}$, the normal equation of (\ref{Tikhonov}) is
\begin{equation}
\left(\A^{\T}*\A+\mu^{-1}\LL^{\T}*\LL\right)*\overrightarrow{\X}=\A^{\T}*\overrightarrow{\B}.
\label{Eequation}
\end{equation}
By Theorem 4.1, substituting the T-GSVD of $\A$ and $\LL$
\begin{equation*}
\A=\U*\D_{\A}*\M^{-1},\quad\LL=\V*\D_{\LL}*\M^{-1},
\end{equation*}
into (\ref{Eequation}), the regularized solution $\overrightarrow{\X}_{\mu}$ can be expressed as
\begin{equation}
\begin{split}
\overrightarrow{\X}_{\mu}&=
\sum_{i=1}^{m-2}
\overrightarrow{\M}_i*\left(
\left(\D_{\A}\right)_{i,i}^{\T}*\left(\D_{\A}\right)_{i,i}
+\mu^{-1}\left(\D_{\LL}\right)_{i,i}^{\T}*\left(\D_{\LL}\right)_{i,i}
\right)^{-1}
*\left(\D_{\A}\right)_{i,i}
*\left(\overrightarrow{\U}_i^{\T}*\overrightarrow{\LL}\right)\\
&+\sum_{i=m-1}^{m}\overrightarrow{\M}_i
*\left(\overrightarrow{\U}_i^{\T}*\overrightarrow{\LL}\right).\\
\end{split}
\label{solution}
\end{equation}
This expression is of the similar form as the matrix cases \cite[equation (2.3)]{Wei16} and is important for our numerical examples in section 5.

\section{Numerical examples}

\hskip 2em In this section, we report results of our T-GSVD method for three test problems from \cite{Reichel21}, where Reichel and Ugwu used them for methods global t-product Arnoldi-Tikhonov regularization (tGAT) and generalized global t-product Arnoldi-Tikhonov regularization (G-tGAT). Numerical experiments show that our proposed T-GSVD method can also be used to solve the Tikhonov regularization problem (\ref{Tikhonov}), and the relative error is of almost the same magnitude as that of the tGAT and G-tGAT methods. All computations were carried out in MATLAB 2019a on an Asus computer with Inter Core i5 processor and 8 GB RAM.

\hskip 2em In image restoration problems, $\overrightarrow{\B}_{\true}\in\RR^{m \times 1 \times n}$ is unavailable error-free data tensor and $\overrightarrow{\B}\in\RR^{m \times 1 \times n}$ in (\ref{minF}) is available but contaminated by the noise tensor $\overrightarrow{\E}\in\RR^{m \times 1 \times n}$, $i.e.$
\begin{equation}
\overrightarrow{\B}=\overrightarrow{\B}_{\true}+\overrightarrow{\E}.
\label{B}
\end{equation}

\hskip 2em In our all examples, the entries of $\overrightarrow{\E}$ are normally distributed random numbers with zero mean and are scaled to correspond to a specified noise level $\tilde{\delta}$. Thus,
\begin{equation}
\overrightarrow{\E}=\tilde{\delta}
\frac{\overrightarrow{\E}_0}{\parallel\overrightarrow{\E}_0\parallel_F}
\parallel\overrightarrow{\B}_{\true}\parallel_F,
\label{noise}
\end{equation}
where the entries of the error tensor $\overrightarrow{\E}_0$ obey $N(0,1)$. $\overrightarrow{\X}_{\true}\in\RR^{m \times 1 \times n}$ is the unknown exact solution of minimal Frobenius norm of the unavailable linear system
\begin{equation*}
\A*\overrightarrow{\X}=\overrightarrow{\B}_{\true}.
\end{equation*}

\hskip 2em Suppose $\overrightarrow{\X}_{\mu}$ is given by (\ref{solution}) and the relative error
\begin{equation*}
E=\frac{\parallel\overrightarrow{\X}_{\mu}-\overrightarrow{\X}_{\true}\parallel_F}
{\parallel\overrightarrow{\X}_{\true}\parallel_F}
\end{equation*}
is used to determine the effectiveness. In our examples, the final relative error is the average of 10 experiments. All the results are showed in Table 1.

\begin{ex}
\upshape{This} example gives Tikhonov regularization solution implemented by our T-GSVD method and the G-tGAT method in \cite[example 6.1]{Reichel21} with the regularization tensor $\LL_2\in\RR^{255 \times 256 \times 256}$.  The blurred tensor $\A\in\RR^{256 \times 256 \times 256}$ is generated as follows
\begin{equation*}
\begin{split}
K_1&=\rm{gravity}\left(256,1,0,1,d\right),\quad d=0.8,\\
K_2&=\rm{gallery}\left('prolate',256,\alpha\right),\quad \alpha=0.46,\\
A_i&=K_1\left(i,1\right)K_2,\quad i=1,2,\dotsc,256,
\end{split}
\end{equation*}
\end{ex}
where the function 'gravity' is from Hansen's Regularization Tools \cite{Hansen07} and $K_2$ is a symmetric positive definite ill-conditioned Toeplitz matrix.  The exact data tensor is $\B_{\true}=\A*\X_{\true}$ and $\X_{\true}\in\RR^{256\times 3\times 256}$ has all entries equal to unity. Set noise level $\tilde{\delta}=10^{-3}$. By using the T-GSVD of $\left\{\A , \LL\right\}$ and combining the Algorithm 4.1 and (\ref{solution}), we can get the approximate solution.

\hskip 2em Note that, in this example, the number of lateral slices of $\X_{\true}$ is $3$. Then we can use (\ref{solution}) to every lateral slices of $\X_{\true}$ to get the sloution $\X_{\mu}$. See \cite{Reichel21} for the choice of regularization parameter $\mu$.

\begin{ex}
\upshape{(2D image restoration)} The true image
$\rm{\bf{Telescope}}$ \cite[example 6.2]{Reichel21} of size $300\times 300$ pixels has been contaminated by blur and noise. The blurred tensor $\A\in\RR^{300 \times 300 \times 300}$ is generated with the MATLAB functions as follows
\begin{equation*}
\begin{split}
&Z_1=\left[\rm{exp}\left(-\left(\left[0:\rm{band\ -\ 1}
\right].^2\right)/\left(2\sigma^2\right)\right),\rm{zeros}\left(1,N\ -\ \rm{band}\right)\right],\\
&K_2=\frac{1}{\sigma\sqrt{2\pi}}\rm{toeplitz}\left(Z_1\right),\\
&Z_2=\left[Z_1\left(1\right)\ \rm{fliplr}\left(Z_1\left(\rm{end}-\rm{length}\left(Z_1\right)+2:\rm{end}\right)\right)\right],\\ &K_1=\frac{1}{\sigma\sqrt{2\pi}}\rm{toeplitz}\left(Z_1,Z_2\right),\\
&A_i=K_1\left(i,1\right)K_2,\quad i=1,2,\dotsc,300,\\
\end{split}
\end{equation*}
where $N=300$, $\sigma=3$ and $\rm{band}=9$. We use (\ref{B}) and (\ref{noise}) to get the blurred and noisy image $\overrightarrow{\B}\in\RR^{300\times 1\times 300}$. Then, the restored image can be get by T-GSVD method and tGAT method \cite{Reichel21} with relative error in Table 1.
\end{ex}

\begin{ex}
\upshape{(Color image restoration)} The true color image $\rm{\bf{flower}}$ \cite[example 6.3]{Reichel21} of size $300\times 300$ pixels has been contaminated by blurred tensor $\A\in\RR^{300 \times 300 \times 300}$ which is the same as in Example 5.2 but $\rm{band}=12$.
The true image is also contaminated by noisy tensor generated by (\ref{B}) and (\ref{noise}). Noise level $\tilde{\delta}=10^{-3}$ is fixed. The restored image by our T-GSVD method is almost the same as that by G-tGAT method after 32 iterations.
\end{ex}

\vspace{-0.3cm}
\begin{table}[htb]
\caption{Results for Examples}
\centering  
\begin{tabular}{ccccc}
\toprule
                          &Noise Level                  &$\mu$
                   &Method                 &Relative error                   \\\toprule
\multirow{2}*{Ex. 5.1}    &\multirow{2}*{$10^{-3}$}     &\multirow{2}*{$7.13e-02$}
                   &G-tGAT                 &$0.00620$                        \\
~                         &~                            &~
                   &T-GSVD                 &$0.01841$                        \\
\specialrule{0em}{0pt}{4pt}
\multirow{2}*{Ex. 5.2}   &\multirow{2}*{$10^{-3}$}     &\multirow{2}*{$3.18e+04$}
                   &tGAT                   &$0.13400$                           \\
~                         &~                            &~
                   &T-GSVD                 &$0.12649$                           \\
\specialrule{0em}{0pt}{4pt}
\multirow{2}*{Ex. 5.3}   &\multirow{2}*{$10^{-3}$}     &\multirow{2}*{$7.34e+03$}
                   &G-tGAT                 &$0.0667$                          \\
~                         &~                            &~
                   &T-GSVD                 &$0.0671$                          \\
\toprule
\end{tabular}
\end{table}
\vspace{-0.3cm}

\hskip 2em The method t-GAT works with the lateral slices of the data tensor $\B$ independently while G-tGAT works with those simultaneously. They are all based on global t-Arnoldi process \cite{Reichel21}, and the first step for the two methods is to reduce the tensor $\A$ to upper Hessenberg matrix. The dimension of the upper Hessenberg matrix must be determined through iterating. Different from tGAT and G-tGAT, our method T-GSVD belongs to direct method. As long as the decomposition structure is obtained by Algorithm 4.1, we can directly compute the solution by (\ref{solution}). The relative errors of different methods are of almost the same magnitude. Hence, our T-GSVD provides a new approach for solving tensor Tikhonov regularization problem. But, T-GSVD may be costly when the data is large-scale. For further improvements of computational complexity, truncated T-GSVD or randomized T-GSVD will be studied in later work.

\section{Conclusion}

\hskip 2em This paper extends the CS decomposition for orthogonal matrices to orthogonal tensors and the GSVD for two matrices to two tensors. Applying the T-GSVD to the data tensor $\A$ and regularization operator $\LL$ in the Tikhonov regularization problem, we have derived the direct expression for the unique solution. Finally, numerical examples illustrate that our method is effective while solving image restoration problems. The restored image determined by T-GSVD is at least as good as other methods, such as tGAT and G-tGAT.

\hskip 2em For matrix cases, the CS decomposition can be used to prove the GSVD. Naturally, how to prove the T-GSVD by the T-CSD deserves more attention. The key to solve this problem is to deal with the tubal scalars which are nonzero but not invertible. When the scale of the data tensor is large in applications, in order to accelerate speed and save memory, we may consider techniques such as truncation and randomization, which will be a future research topic.


\end{document}